\documentclass[reqno,12pt]{amsart}
\usepackage{latexsym,amssymb,amsthm,amsmath}

\theoremstyle{plain}

\theoremstyle{remark}
\newtheorem{remark}{\bf Remark}

\numberwithin{equation}{section}

\def\<{\left < }
\def\>{\right >}
\def\({\left ( }
\def\){\right )}

\begin{document}

\title[OPEN PROBLEMS AND CONJECTURES Revisited]{Some open problems and conjectures on submanifolds of finite type: revisited}

\author[B.-Y. Chen]{Bang-yen Chen}
\vskip.15in
\begin{abstract} Submanifolds of finite type were introduced by the author during the late 1970s (cf. [C1-C4]). The first results on this subject had been collected in author's books
[C4,C7]. A detailed survery on results on finite submanifolds up to 1996 was given in \cite{c3}. Moreover, a list of ten open problems and three conjectures on submanifolds of finite type was published in  \cite[1981]{C18}. 
The main purpose of this article is thus to provide some updated information on the three conjectures listed in \cite{C18}. \end{abstract} 

\address{Department of Mathematics, \\
	Michigan State University \\East Lansing, Michigan 48824--1027, U.S.A.}
\email{bychen@math.msu.edu}

\maketitle

\vskip.2in
\section{Preliminaries}

 Let $x:M\rightarrow \mathbb E^m$ be an immersion of an $n$-dimensional, connected manifold $M$ into the Euclidean $m$-space $\mathbb E^m$. With respect to the Riemannian metric $g$ on $M$ induced from the Euclidean metric  of the ambient space $\mathbb E^m$, $M$ is a
Riemannian manifold. Denote by $\Delta$ the Laplacian operator of the Riemannian manifold $(M,g)$. The immersion $x$ is said to be {\it of finite type\/}  if each component of the position vector field of $M$ in $\mathbb E^m$, also denoted by $x$, can be written as a finite sum of eigenfunctions of the Laplacian operator, {\it that is,\/} if
\begin{align}x=c+x_{1}+x_{2}+\ldots+x_{k}\end{align}
where $c$ is a constant vector, $x_{1},\ldots,x_{k}$ are non-constant maps satisfying $\Delta x_{i}=\lambda_{i}x_{i}, \; i=1,\ldots,k.$ In particular, if all eigenvalues
$\{\lambda_{1},\ldots,\lambda_{k}\}$ are mutually different, then the immersion $x$ (or the submanifold $M$) is said to be of {\it k-type\/} . And the decomposition (1.1) is called the {\it spectral decomposition\/} (or the {\it spectral resolution}) of the immersion $x$. In particular, if one of  $\{\lambda_{1},\ldots,\lambda_{k}\}$ is zero, then the immersion is said to be of {\it null k-type.}  Just like  minimal submanifolds, submanifolds of finite type can be described by a spectral variation principle, namely as critical points of directional deformations (see, \cite{CDVV2,CDVV4,book} for details).

A submanifold is said to be {\it of infinite type\/} if it is not of finite type. It is clear that every submanifold of null $k$-type is non-compact. If $x(M)$ is contained in a hypersphere $S^{m-1}$ of $\mathbb E^m$, then the immersion $x$ is said to be {\it mass-symmetric\/} in $S^{m-1}$ if the constant vector $c$ in the spectral decomposition (1.1) is the center of the hypersphere $S^{m-1}$ in $\mathbb E^m$. 

In terms of finite type submanifolds, a  result of Takahashi [Ta1] says that a submanifold of a Euclidean $m$-space $\mathbb E^m$ is of 1-type if and only if it is either a minimal
submanifold of $\mathbb E^m$ or it is a minimal submanifold of a hypersphere of the Euclidean space.

 Let $M$ be a finite type submanifold whose spectral
decomposition is given by (1.1). If we define a polynomial
$P$ by
\begin{align}P(t)=\prod_{i=1}^{k} (t-\lambda_{i}),\end{align}  then $P(\Delta)(x-c)=0.$ 
The polynomial $P$ is called the {\it minimal
polynomial\/} of the finite type submanifold $M$. For
an $n$-dimensional submanifold $M$ of a Euclidean space,
the mean curvature vector $H$ of $M$ satisfies 
\begin{align}\Delta x = -nH.\end{align}
From (1.3) we see that the minimal polynomial $P$
of the k-type submanifold $M$ 
also satisfies the condition $P(\Delta)H=0$. 
Conversely,  if $M$ is compact and if there
exists a constant vector $c$ and a nontrivial polynomial
$P$ such that $P(\Delta)(x-c)=0$ (or $P(\Delta)H=0$), then
$M$ is of finite type (see, \cite{C2,C4,C7}.)  This
characterization of finite type submanifolds in terms of
the minimal polynomial plays an important role in the
study  of submanifolds of finite type. 

If $M$ is non-compact, then the existence of a nontrivial polynomial
$P$ such that that $P(\Delta)H=0$ does not imply that $M$
is of finite type in general. However, if $M$ is
1-dimensional or $P$ is a polynomial  of degree $k$ which
has exactly $k$ distinct  roots,  then the existence
of the polynomial $P$ satisfying the condition
$P(\Delta)(x-c)=0$ for constant vector $c$ guarantees that
$M$ is of finite type (in fact, it is of $k$-type with $k
\leq deg(P)$) (see, \cite{CP}).

The following formula of $\Delta H$ obtained in \cite{C2,C4,C7} plays an important role in the study of
submanifolds of low type  and also in the study of biharmonic submanifolds (cf. \S 5.1). 
\begin{align}\Delta
H=\Delta^{D}H + \sum_{i=1}^{n} h(e_{i},A_{H}e_{i}) + 2\,tr(A_{DH})+ {n\over 2} grad \<H,H\> ,\end{align} 
where
$\Delta^{D}$ is the Laplacian operator associated with the
normal connection $D$, $h$ the second fundamental form,
and $\{ e_{1},\ldots,e_{n}\}$ a local orthonormal frame of
$M$. In particular, if $M$ is a hypersurface of a
Euclidean space $\mathbb E^{n+1}$ (respectively, if $M$ is a
hypersurface of the unit hypersphere of $\mathbb E^{n+2}$ centered
at the origin), then formula (1.4) yields 
\begin{align}\Delta H=(\Delta \alpha +\alpha ||h||^{2})\xi + 2\, tr(A_{DH})+{n\over 2}\, grad\<H,H\>,\end{align}
\begin{align} & (respectively,\;\;\;\Delta H=(\Delta \alpha' +\alpha'
||h||^{2})\xi - n\<H,H\>x \\ &\notag\hskip1.5in
+2\,tr(A_{DH})+{n\over 2}\, grad \<H,H\>,)\end{align}
where $\alpha$ is the mean curvature and $\xi$ a unit
normal vector of $M$ in $\mathbb E^{n+1}$ (respectively, $\alpha'$
is the mean curvature and $\xi$ a unit normal vector
 of $M$ in $S^{n+1}$.)  

A  similar formula of $\Delta H$ holds for submanifolds in pseudo-Euclidean spaces. Also, similar formulas of $\Delta H$ for submanifolds in the hyperbolic hypersurface $H^{m-1}$ and in the anti-de Sitter hypersurface $H^{m-1}_{1}$ also had been obtained in \cite{C6,C9}.

\vskip.4in 
\section{Submanifolds of Finite Type} 
\subsection {Finite Type Hypersurfaces of Euclidean Space.}

The class of finite type submanifolds is very large. For
instances, minimal submanifolds of a Euclidean space,
minimal submanifolds of a hypersphere are of 1-type and
compact homogeneous submanifolds, equivariantly immersed,
are of finite type \cite{C4} (see, also [De1,De2] and [Ta2]
for irreducible case).  But very few is know about the most
elementary submanifolds of Euclidean space, namely
hypersurfaces of a Euclidean space, in particular,
surfaces in Euclidean 3-space. So far no surfaces of
finite type in $\mathbb E^3$ are known, other than minimal
surfaces,  circular cylinders and the spheres. Therefore,
the following problem seems to be very interesting (cf.
\cite{C7,C11}). 
\vskip.1in

 {\bf Problem 1.} {\it Classify all finite type
hypersurfaces in $\mathbb E^{n+1}$. In particular, classify all
finite type surfaces in $\mathbb E^3$.}
\vskip.1in

Related to this problem we recall that every compact
2-type hypersurface of $\mathbb E^{n+1}$ has non-constant mean
curvature  \cite{CLu} (this had also been pointed out by
Garay, see \cite{CLu}). And if $n=1$, this problem has a
complete solution. In fact, it was proved in \cite{C3,C4}
that circles are the only finite type closed planar
curves.  It was pointed out in \cite{C14} that lines are the
only non-closed planar curves of finite type (in fact, it
is the only null finite type planar curves) (cf. \cite{CDVV}
for the details.)  

The first result concerning the
classification of finite type surfaces in Euclidean
3-space was obtained in \cite{C11}  which stated
 that  circular cylinders are the only tubes  of
finite type. 
In \cite{G2} it shown that a cone in $\mathbb E^3$ is of finite type if
and only if it is a plane. It was proved  in \cite{CDVV} that a
ruled surface in $\mathbb E^3$ is of finite type if and only if it
is a plane, a circular cylinder or a helicoid. Furthermore,
\cite{CD2} proved that  spheres and  circular cylinders are the
only quadrics of finite type in $\mathbb E^3$. It seems
 that the only surfaces of finite type in
$\mathbb E^3$ are minimal surfaces, spheres and circular
cylinders. For compact finite type surfaces in
$\mathbb E^3$, the author made the following conjecture in
\cite{C7,C18}.
\vskip.1in

 {\bf Conjecture 1.} {\it The only compact
finite type surfaces in $\mathbb E^3$ are the spheres.}
\vskip.1in

Some partial affirmative answers to this conjecture were given in \cite{B,AGM,C11,DDV1,DPV,HV4,V2}. However, Conjecture 1 remains open.

It also seems to the author that hyperspheres are the only
compact  hypersurfaces of finite type in Euclidean space.
However, this seems to be a quite difficult problem.

\subsection{Finite Type Hypersurfaces of Hyperspheres.}

For finite type hypersurfaces of a hypersphere the
following problem is also quite interesting.
\vskip.1in

{\bf Problem 2.} {\it Classify finite type
hypersurfaces of a hypersphere in $\mathbb E^{n+2}$.}
\vskip.1in

In contrast with hypersurfaces of finite type in Euclidean
space, there exist many examples of 1-type hypersurfaces as well as many examples of 
(mass-symmetric) 2-type hypersurfaces in a hypersphere of
$\mathbb E^{n+2}$. In fact, it was first proved in \cite{C8} that every
isoparametric hypersurface of a hypersphere is either of
1-type or mass-symmetric and of 2-type. Since there are
ample  examples of non-minimal isoparametric hypersurfaces
in a hypersphere, we have ample examples of mass-symmetric
2-type hypersurfaces  in a hypersphere.  It was proved
that every 3-type spherical hypersurface  has non-constant
mean curvature (\cite{C15}, \cite{CL}).   Moreover, 
\cite{CD1} proved that standard 2-spheres in $S^3$ and  products of  plane circles are the only finite type compact surfaces with constant Gauss curvature in $S^3$.
From all the information we have, it seems to the author
that there exist no surfaces of k-type in $S^3$ for any finite $k$ greater than 2. 

 At an
international conferences held at Berlin in  the summer of
1990 the author had 
  announced the following conjecture concerning finite type
surface in a 3-sphere.
\vskip.1in

 {\bf Conjecture 2.} {\it Minimal surfaces,
standard 2-spheres and products of plane circles are the
only finite type surfaces in $S^3$} (imbedded
standardly in $\mathbb E^4$).
\vskip.1in

Some partial affirmative answers to this conjecture were obtained in \cite{C15,CBG,CL,FL,Na}.
\begin{remark}This conjecture stays open.\end{remark}

\vskip.0in
\section{2-type Submanifolds}
\vskip.05in 

2-type submanifolds are the ``simplest submanifolds'' next to minimal submanifolds. In particular, 2-type submanifolds mass-symmetric spherical 2-type submanifolds, deserve  special attention.  

Mass-symmetric spherical 2-type submanifolds have some special properties. For instances, every mass-symmetric spherical 2-type submanifolds has constant mean curvature which is completely determined by its order \cite{C4}. Moreover,
such a submanifold is  {\it pointwise orthogonal\/} (in the sense of \S 4.)

\subsection{Spherical 2-type hypersurfaces.}
The first  classification theorem of 2-type surfaces of a hypersphere was obtained in \cite{C3,C4} which says that a compact surface  of a hypersphere $S^3$ in
$\mathbb E^4$ is the product of two plane circles with different radii if and only if it is mass-symmetric and of 2-type. \cite{BG1} proved that the  the assumption on ``mass-symmetry'' can be removed. It was proved in  \cite{HV1}  that the above classification of 2-type
surfaces in $S^3$   indeed is of local nature which says that the same result still holds without the assumption of compactness.  

For a 2-type hypersurface $M$ of a hypersphere $S^{n+1}$,  it was proved in \cite{C8} that
every mass-symmetric 2-type hypersurface $M$ of $S^{n+1}$ has nonzero constant mean curvature in $S^{n+1}$ and constant scalar curvature (this result is also of
local nature,  because the proof of this fact given in \cite{C8} did not  use compactness of $M$  at all). 

It was shown in \cite{CBG} that a spherical 2-type hypersurface is mass-symmetric if and only if it has constant mean curvature. Locally, a hypersurface $M$ of
$S^{n+1}$ is the product of two spheres $S^{p}(r_{1})\times S^{n-p}(r_{2})$ with $r_{1}^{2}+r_{2}^{2}=1$ and $(r_{1},r_{2})\not= (\sqrt{p\over n},\sqrt{{n-p}\over n})$
if and only  if $M$ is of 2-type and it has at most two distinct principal curvatures \cite{C16} (see, also \cite{CBG} for compact case). Recently,  \cite{HV2} proved that  every 2-type
hypersurface of a hypersphere $S^{n+1}$ is mass-symmetric  and hence it has constant mean curvature and constant scalar curvature .

Conversely, it was known that each hypersurface of $S^{n+1}$ with nonzero constant mean curvature in $S^{n+1}$ and constant scalar curvature is mass-symmetric and of 2-type unless it is a portion of small hypersphere of $S^{n+1}$ (a result obtained in \cite{C8} for compact case. For noncompact case, this fact is a consequence of  formula (1.6) and Proposition 4.3 of \cite{CP} (see  \cite{CL}) and this  was also pointed out  in \cite{HV2} independently). For spherical 2-type hypersurfaces, the author had proposed in \cite{C18} the following 
\vskip.1in

 {\bf Problem 3.} {\it Study and classify 2-type hypersurfaces in a hypersphere of $\mathbb E^{n+2}$. In particular, classify 3-dimensional 2-type hypersurfaces of a hypersphere $S^4$ in $\mathbb E^5$.}
\vskip.1in

Since every spherical 2-type hypersurface is mass-symmetric, a result of \cite{C8} implies that if either $M$ is a 3-dimensional 2-type Dupin hypersurface of
$S^4$ or  $M$ is a spherical 2-type Dupin hypersurface
with at most 3 distinct principal curvatures, then the
Dupin hypersurface is isoparametric.  For a general Dupin hypersurfaces,
the author asked in \cite{C18} the following
\vskip.1in

{\bf Problem 4.} {\it When is a Dupin hypersurface $M$ of  a hypersphere  of finite type? When is a finite type Dupin hypersurface of a hypersphere  isoparametric?}
\vskip.1in

 In contrast with the existence of many compact 2-type hypersurfaces in hyperspheres of a Euclidean space,    {\it there exits no compact 2-type hypersurface in the hyperbolic space $H^{n+1}$} (imbedded standardly in the Minkowski space-time $\mathbb E^{n+2}_1$
 by the equation $\<x,x\>=-1,\,\, t>0$), although there exist complete, non-compact 2-type hypersurfaces in  $H^{n+1}$. (See  \cite{C9,C19} for more results on finite type submanifolds in pseudo-Euclidean spaces, in particular, on finite type submanifolds in
 hyperbolic spaces and in de Sitter space-times. For instances, 2-type surfaces in hyperbolic 3-spaces have been classified in \cite{C19} and it is also known in \cite{C19}
that {\it every space-like 2-type hypersurface of the de Sitter space-time $S^{n+1}_{1}$, imbedded standardly in $\mathbb E^{n+2}_{1}$, has nonzero constant mean curvature and
constant scalar curvature.\/}) 

\vskip.1in

\subsection{2-type Submanifolds of Codimension 2.}

 It was proved in \cite{BC1} that there  exist no  mass-symmetric surfaces in $S^4$. So far there are no known examples of non-mass-symmetric 2-type surfaces in $S^4$. In this respect,  the author   asked in \cite{C18} the following 
\vskip.1in

{\bf Problem 5.} {\it Do there exist (non-mass-symmetric) 2-type surfaces in $S^4$?} 
\vskip.1in

 \cite{CL} showed that  a  2-type submanifold  with parallel mean curvature vector
is either spherical  or of null 2-type. Related with this fact, the author had
 asked in \cite{C18} the following 
\vskip.1in

{\bf Problem 6.} {\it Is every $n$-dimensional  non-null 2-type submanifold of $\mathbb E^{n+2}$   with constant mean curvature  spherical?}
\vskip.1in

If the answer to Problem 6 is affirmative, then an $n$-dimensional 2-type submanifold of $\mathbb E^{n+2}$ is spherical if and only if it is non-null and it has constant mean curvature.

\vskip.1in

\subsection{2-type Spherical Surfaces of
Higher Codimension.}

There exist ample examples of 1-type (mass-symmetric) surfaces which lies fully in
odd-dimensional spheres as well as in even-dimensional spheres.  Although there exists abundant examples of mass-symmetric 2-type surfaces lying fully in  odd-dimensional hyperspheres (cf. for instances, \cite{BC1,C4, G1, K, Mi}), in contrast there exist no  examples of mass-symmetric 2-type surfaces which lie fully in a  hypersphere of a
Euclidean space for any even codimension. Hence, it is natural to ask the following problem  which is more general than Problem 5 (see, also \cite{K}).
\vskip.1in

{\bf Problem 7.} {\it Do there exist  2-type surfaces which lie fully in an even-dimensional
hypersphere of a Euclidean space? In particular, do there exist mass-symmetric 2-type surfaces which lie fully in an even-dimensional hypersphere?}
\vskip.1in

The author would like to point out
that the answer to Problem 7 is negative if the
mass-symmetric spherical 2-type surface is one of the
following surfaces: 

\begin{itemize}
\item a stationary surface \cite{BC1}, 

\item a topological 2-sphere [Ko1], 

\item a surface with constant Gauss curvature \cite{Mi},  

\item a flat Chen surface \cite{G1}.
\end{itemize}
 \vskip.2in

\section{Linearly Independent  submanifolds.}

The notion of linearly independent immersions and
linearly independent submanifolds 
were defined as follows \cite{C17}:

 Let $x:M\rightarrow \mathbb E^m$ be a $k$-type isometric
immersion  whose  spectral decomposition is given by
(1.1).  Denote by $E_{i}$ 
the subspace of $\mathbb E^m$ spanned by $\{ x_{i}(p),\, p \in
M\} \,\,(i\in \{1,\ldots,k\} )$.  The immersion $x$ (or the submanifold $M$) is said
to be {\it linearly independent\/} if the subspaces
$E_{1},\ldots,E_{k}$ are linearly independent. And the
immersion $x$ (or the submanifold) is said to be {\it
orthogonal} if the subspaces $E_{1},\ldots,E_{k}$ are
mutually orthogonal. 

Clearly, every orthogonal immersion
is a linearly independent immersion and every 1-type
immersion is an orthogonal immersion.  There exist
many examples of  orthogonal immersions and abundant
examples of linearly independent immersions which are not
orthogonal. In fact every k-type curve lying fully in $\mathbb E^{2k}$ and
every null $k$-type curve lying fully in $\mathbb E^{2k-1}$ are
linearly independent curves, but  $W$-curves are the only
orthogonal curves in a Euclidean space. 

  For a linearly independent immersion $x:M\rightarrow \mathbb E^m$ and a point $p\in M$ one has the notion of the {\it adjoint hyperquadric\/} $Q_p$ at $p$ (cf. \cite{C17}.) If 
submanifold $M$ lies in one of the adjoint hyperquadrics $Q_p,$  $(p\in M)$, then all of the adjoint hyperquadrics $Q_p,\, (p\in M)$ are the same adjoint hyperquadric. This common adjoint hyperquadric  is  called the {\it adjoint hyperquadric of the linearly independent immersion $x$.}

It was shown in \cite{C17} that if  $x: M\rightarrow \mathbb E^m$ is 
a linearly independent immersion 
of a compact manifold $M$, then the submanifold $M$ is
contained in its adjoint hyperquadric  if and only if
the submanifold $M$ is spherical
 (with an appropriate center). Furthermore,    a
non-minimal, linearly independent immersion $x: M
\rightarrow \mathbb E^m$ is 
orthogonal  if and only if  $M$ is immersed by $x$ as a
minimal submanifold of the adjoint hyperquadric  \cite{C17}
(for some special orthogonal immersions of compact
manifolds, see also \cite{HV3}). As a consequence it follows
that every orthogonal immersion of a compact manifold  is 
spherical. Moreover, one may also prove that every  compact
homogeneous submanifold, equivariantly immersed in $\mathbb E^m$, 
is orthogonal and hence it is immersed as a minimal
submanifold in its adjoint hyperquadric ( \cite{C17}). 
Linearly independent hypersurfaces of a Euclidean space
are hyperspheres, minimal hypersurfaces or spherical
hypercylinders (\cite{CP, HV4, CDVV2}). (See also
\cite{G3} and \cite{DPV} for some special cases.)
By applying  the   classification theorem of 2-type curves
in Euclidean space obtained in
 \cite{CDV}, we may  conclude that the
only linearly independent curves of codimension 2 in a
Euclidean space are circles, lines and circular helices.
As we already know, there are abundant examples of
linearly independent curves of codimension 3 in Euclidean space.

In views of these, the author had proposed the following two problems concerning linearly
independent immersions.
\vskip.1in

{\bf Problem 8.} {\it Study and classify linearly independent  2-type immersions.}
\vskip.1in

{\bf Problem 9.} {\it Study and classify linearly independent  submanifolds of codimension 2.}
\vskip.1in

Let $x: M \rightarrow \mathbb E^m$ be a $k$-type isometric immersion whose spectral decomposition is given by (1.1). The immersion is said to be {\it pointwise linearly
independent\/} (respectively, {\it pointwise orthogonal\/}) if, for each point $p \in M$, the k vectors $x_{1}(p),\ldots, x_{k}(p)$ are linearly independent (respectively, are orthogonal.) 

The class of pointwise linearly independent submanifolds and pointwise orthogonal submanifolds are much wider than the class of linearly independent submanifolds and orthogonal submanifolds. For example,  every mass-symmetric spherical 2-type submanifold is   pointwise orthogonal, although it is not orthogonal in general. (This follows from the definition of mass-symmetric spherical 2-type submanifolds and Theorem
4.1 of p.274 of \cite{C4} (see \cite{Mi}).)

One may study  similar problems for pointwise linearly independent submanifolds and for
pointwise orthogonal submanifolds. However, these  seem to  be quite difficult.

\vskip.2in

\section{Biharmonic submanifolds.} 

Let $x: M \rightarrow \mathbb E^m$ be an isometric immersion. As we mentioned in Preliminaries,  the position vector of $M$ in $\mathbb E^m$ satisfies
\begin{align}\Delta x=-nH. \end{align}
   Formula (5.1) implies that the immersion is minimal if and only if  the immersion 
is harmonic, {\it that is,\/} $\Delta x = 0.$ An isometric immersion $x : M \rightarrow E^m$ is called {\it biharmonic\/}  if we have
\begin{align}\Delta^{2}x=0,\hskip.2in that\,\, is,  \hskip.2in \Delta H=0.\end{align}
It is obvious that  minimal immersions are biharmonic.

 In \cite{C18}  the author had asked the following simple geometric question (see also [CI1]). 
\vskip.1in

{\bf Problem 10.} {\it Other than minimal submanifolds of
$\mathbb E^m$, which submanifolds of $\mathbb E^m$ are biharmonic?}
\vskip.1in

The study of biharmonic submanifolds was initiated by the author in the middle of 1980s in his program of understanding the finite type submanifolds in Euclidean spaces; also independently by G.-Y. Jiang  \cite{J} for his study of Euler-Lagrange's equation of bienergy functional in the sense of Eells and Lemaire. 

The author showed in 1985 that biharmonic surfaces in $\mathbb E^3$ are minimal (unpublished then, also independently by Jiang \cite{J}). This result was the starting point of I. Dimitric's work on his doctoral thesis at Michigan State University (cf. \cite{Dim89}).  
In particular, Dimitric extended author's unpublished result to  biharmonic hypersurfaces of $\mathbb E^m$ with at most  two distinct principal curvatures \cite{Dim89}.  In his thesis,
Dimitric also proved that every biharmonic submanifold of finite type 
in ${\mathbb{E}}^m$ is minimal.
Another extension of this result on biharmonic surfaces was given by T. Hasanis and T. Vlachos in \cite{HV} (see also
 \cite{D98}). They proved that biharmonic hypersurfaces of ${\mathbb{E}}^4$ are minimal.

The author made in \cite{C18} the following Biharmonic Conjecture.
\vskip.06in

\hskip.3in {\bf Conjecture 3}:
\emph{The only biharmonic submanifolds of Euclidean spaces are 
the minimal ones.}
\vskip.06in

A {\it biharmonic map} is a map $\phi:(M,g)\to (N,h)$ between Riemannian manifolds that is a critical point of the bienergy functional:
\begin{align}E^2(\phi,D)=\frac{1}{2}\int_D ||\tau_\phi||^2* 1\end{align}
for every compact subset $D$ of $M$, where
$\tau_\phi={\rm trace}_g\nabla d\phi$ is the tension field $\phi$. 
  The Euler-Lagrange equation of this functional gives the biharmonic map
equation (see \cite{J})
\begin{align}\label{BH} \tau^2_\phi:={\rm trace}_g(\nabla^\phi\nabla^\phi-\nabla^\phi_{\nabla^M})\tau_\phi-{\rm trace}_g R^N(d\phi,\tau_\phi)d\phi=0,\end{align}
where $R^N$ is the curvature tensor of $(N,h)$.
Equation \eqref{BH} states that $\phi$ is a biharmonic map if and only if its bi-tension field $\tau^2_\phi$ vanishes.  

Let $M$ be an $n$-dimensional submanifold of a Euclidean $m$-space $\mathbb E^m$.
 If we denote by $\iota:M\to \mathbb E^m$ the inclusion map of the submanifold, then the tension field of the inclusion map is given by $\tau_\iota=-\Delta\iota=-n H$ according to Beltrami's formula. Thus $M$ is a biharmonic submanifold if and only if 
$n\Delta  H=- \Delta^2 \iota=-\tau^2_\iota=0,$
i.e., the inclusion map $\iota$ is a biharmonic map.

 Caddeo,  Montaldo and  Oniciuc \cite{CMO02} proved that every biharmonic surface in 
the hyperbolic $3$-space $H^3(-1)$ of constant curvature $-1$ is minimal. They also proved that biharmonic hypersurfaces of
$H^n(-1)$ with at most two distinct principal curvatures are minimal \cite{CMO01}. 
Based on these, Caddeo, Montaldo and Oniciuc made in  \cite{CMO01} the following {\it generalized biharmonic conjecture.}
\vskip.06in

\hskip.3in  
\emph{Every biharmonic submanifold of a Riemannian manifold with non-positive sectional curvature is minimal.}
\vskip.06in

The study of  biharmonic submanifolds is nowadays a very active research subject. In particular, since 2000 biharmonic submanifolds have been
receiving a growing attention and have become a popular subject of study with
many progresses.

\subsection{Recent developments  on my original biharmonic conjecture}

 Let  $x :M\to \mathbb E^m$ be an isometric immersion of a Riemannian $n$-manifold $M$ into a Euclidean $m$-space $\mathbb E^m$.  
 Then $M$ is biharmonic if and only if it satisfies the following fourth order strongly elliptic semi-linear PDE system (see, for instance, \cite{C2,C3,C4,book})
 \begin{align}\notag &\begin{cases}\; \Delta^{D} H + \sum_{i=1}^{n}
\sigma(A_{H}e_{i},e_{i})=0,\\  \\ \;n\,\nabla\!\left<\right.\! H, H\! \left.\right> + 4\, {\rm trace}\,A_{D H}=0,\end{cases}\end{align}
where $\Delta^D$ is the Laplace operator associated with the normal connection $D$, $\sigma$ the second fundamental form, $A$ the shape operator, $\nabla\!\left<\right.\! H, H\! \left.\right>$ the gradient of the squared mean curvature, and $\{e_1,\ldots,e_n\}$ an orthonormal frame of $M$.

An immersed submanifold $M$ in a Riemannian manifold $N$ is said to be {\it properly
immersed} if the immersion is a proper map, i.e.,  the preimage of each compact set in
$N$ is compact in $M$.

The {\it total mean curvature} of a submanifold $M$ in a Riemannian manifold is given by $\int_M | H|^2 dv$.

 Denote by $K(\pi)$ the sectional curvature of a given Riemannian $n$-manifold $M$ associated with a plane section $\pi\subset T_pM$, $p\in M$. For any orthonormal basis $e_1,\ldots,e_n$ of the tangent space $T_pM$, the scalar curvature $\tau$ at $p$ is defined to be $\tau(p)=\sum_{i<j} K(e_i\wedge e_j). $

Let $L$ be a subspace of $T_pM$  of dimension $r\geq 2$  and $\{e_1,\ldots,e_r\}$ an orthonormal basis of $L$. The scalar curvature $\tau(L)$ of $L$ is defined by
$$\tau(L)=\sum_{\alpha<\beta} K(e_\alpha\wedge e_\beta),\quad 1\leq \alpha,\beta\leq r.$$

For an integer  $r\in [2,n-1]$, the {\it $\delta$-invariant} $\delta(r)$ of $M$ is defined  by  (cf. \cite{C20,C00,book})
\begin{align}\label{1.3} \delta(r)(p)=\tau(p)- \inf\{\tau(L)\},\end{align} where $L$ run over all $r$-dimensional linear subspaces of $T_pM$. 

For any $n$-dimensional submanifold $M$ in $\mathbb E^m$ and any integer $r\in [2, n-1]$,  the author proved the following general sharp inequality  (cf. \cite{C00,book}):
\begin{align}\label{1.4} \delta(r) \leq  \frac{n^2(n-r)}{2(n-r+1)} |\overrightarrow H|^2.\end{align}

 A submanifold  in $\mathbb E^m$ is called {\it $\delta(r)$-ideal} if it satisfies the equality case of \eqref{1.4} identically. Roughly speaking ideal submanifolds are submanifolds which receive the least possible tension from its ambient space (cf. \cite{C00,C11}).
 
 A hypersurface of a Euclidean space is called {\it weakly convex} if it has non-negative principle curvatures.
 
It follows immediately from the definition of biharmonic submanifolds and Hopf's lemma that every biharmonic submanifold in a Euclidean space is non-compact.

\vskip.06in
The following provides an overview of some affirmative partial solutions to my original biharmonic conjecture.

\begin{itemize}
\item Biharmonic surfaces in $\mathbb E^3$ (B.-Y. Chen \cite{C18,book} and G. Y. Jiang \cite{J}).

\item Biharmonic  curves (I. Dimitric \cite{Dim89,Di1}).

\item Biharmonic  hypersurfaces in $\mathbb E^4$ (T. Hasanis and T. Vlachosin \cite{HV}) (a different proof by F. Defever  \cite{D98}).

\item Spherical submanifolds (B.-Y. Chen \cite{C18}).

\item Biharmonic hypersurfaces with at most 2 distinct principle curvatures (I. Dimitric  \cite{Dim89}).

\item Biharmonic  submanifolds of finite type (I. Dimitric \cite{Dim89,Di1}).

\item Pseudo-umbilical biharmonic  submanifolds (I. Dimitric \cite{Di1}).

\item Biharmonic  submanifolds which are complete and proper (Akutagawa and Maeta  \cite{AM}).

\item Biharmonic properly immersed submanifolds (S. Maeta \cite{M12a}).

\item Biharmonic  submanifolds  satisfying the
decay condition at infinity $$\lim_{\rho\to \infty}\frac{1}{\rho^2}\int_{f^{-1}(B_\rho)}| H |^{2}dv=0,$$
where $f$ is the immersion, $B_\rho$ is a geodesic ball of $N$ with radius $\rho$  (G. Wheeler \cite{Wh}).

\item Submanifolds whose $L^p,\, p\geq 2$,  integral of the mean curvature vector field satisfies
certain decay condition at infinity (Y. Luo  \cite{Luo3}).

\item $\delta(2)$-ideal and $\delta(3)$-ideal biharmonic  hypersurfaces (B.-Y. Chen and M. I. Munteanu  \cite{CM}).

\item Weakly convex biharmonic submanifolds (Y. Luo in \cite{Luo1}).
\end{itemize}

In \cite{Ou09}, Y.-L. Ou constructed examples to show that my original biharmonic conjecture cannot be generalized to the case of biharmonic conformal submanifolds in Euclidean spaces. 
\vskip.05in

\begin{remark}{\bf Conjecture 3 remains open.}\end{remark}

\begin{remark}  Conjecture 3 is false if the ambient Euclidean space were replaced by a pseudo-Euclidean space. The simplest examples are constructed by Chen and Ishikawa in \cite{CI91}. For instance, we have the following.
\vskip.05in

{\bf Example.}  Let $f(u,v)$ be a proper biharmonic function, i.e. $\Delta f\ne 0$ and $\Delta^2 f=0$. Then 
\begin{align}\label{Bi} x(u,v)=(f(u,v),f(u,v),u,v)\end{align}
defines a biharmonic, marginally trapped surface in the Minkowski 4-space $\mathbb E^4_1$ with the Lorentzian metric $g_0=-dt_1^2+dx_1^2+dx_2^2+dx_3^2$. 

Here, by a marginally trapped surface, we mean a space-like surface in $\mathbb E^4_1$ with light-like mean curvature vector field.

It was proved in \cite{CI91} that the biharmonic surfaces defined by \eqref{Bi} are the only biharmonic, marginally trapped surfaces in $\mathbb E^4_1$.
\end{remark}

\subsection{Recent developments  on generalized biharmonic conjecture}

Let $M$ be a submanifold of a Riemannian manifold with inner product $\left<\;\,,\;\right>$, then $M$ is called {\it $\epsilon$-superbiharmonic} if
$$\left<\right.\!\Delta  H, H\! \left.\right>\geq (\epsilon-1)|\nabla H|^2,$$
where $\epsilon\in [0,1]$ is a constant.
For a complete Riemannian manifold $(N,h)$ and $\alpha\geq 0$, if the sectional curvature $K^N$ of $N$ satisfies
$$K^N\geq -L(1+{\rm dist}_N(\,\cdot\,,q_0)^2)^{\frac{\alpha}{2}}$$
for some $L>0$ and $q_0\in N$, 
then we  call that $K^N$ has a polynomial growth bound of order $\alpha$ from below.
\vskip.06in

There are also many affirmative partial answers to the generalized Chen's biharmonic conjecture. The following provides a brief overview of the affirmative partial answers to this generalized conjecture.

\begin{itemize}
\item Biharmonic hypersurfaces in the hyperbolic 3-space $H^3(-1)$ (Caddeo,  Montaldo and  Oniciuc \cite{CMO01}).

\item  Biharmonic hypersurfaces in $H^4(-1)$ (Balmu\c s, Montaldo and  Oniciuc \cite{BMO10b}).

\item Pseudo-umbilical biharmonic submanifolds of $H^m(-1)$ (Caddeo,  Montaldo and  Oniciuc \cite{CMO01}).

\item Biharmonic hypersurfaces of $H^{n+1}(-1)$ with at most two distinct principal curvatures (Balmu\c s, Montaldo and  Oniciuc \cite{BMO08}).

\item Totally umbilical biharmonic hypersurfaces in Einstein spaces (Y.-L. Ou \cite{Ou10}).

\item Biharmonic hypersurfaces with finite total mean curvature in a Riemannian manifold of non-positive Ricci curvature (Nakauchi and Urakawa \cite{NU1}).

\item Biharmonic submanifolds with finite total mean curvature in a Riemannian manifold of non-positive sectional curvature (Nakauchi and Urakawa \cite{NU2}).

\item Complete biharmonic hypersurfaces $M$ in a Riemannian manifold of non-positive Ricci curvature whose mean curvature vector satisfies $\int_M | H|^\alpha dv<\infty$ for some $\epsilon>0$ with $1+\epsilon\leq \alpha<\infty$ (S. Maeta \cite{M13}).

\item Biharmonic properly immersed submanifolds in a complete Riemannian manifold with non-positive sectional curvature whose sectional curvature has polynomial growth bound of order less than 2 from below (S. Maeta \cite{M12b}).

\item Complete biharmonic submanifolds  with finite bi-energy and energy in a non-positively curved Riemannian manifold (N. Nakauchi, H. Urakawa and S. Gudmundsson \cite{NUG}).

\item Complete oriented biharmonic hypersurfaces $M$  whose mean curvature $H$ satisfying $H\in L^2(M)$  in a Riemannian manifold with non-positive Ricci tensor (Al\'{\i}as, Garc\'{\i}a-Mart\'{\i}nez and Rigoli \cite{Al}).

\item Compact biharmonic submanifolds in a Riemannian manifold with non-positive sectional curvature (G.-Y. Jiang \cite{J} and S. Maeta \cite{M13}).

\item $\epsilon$-superbiharmonic submanifolds in a complete Riemannian manifolds satisfying the decay condition at infinity
$$\lim_{\rho\to \infty}\frac{1}{\rho^2}\int_{f^{-1}(B_\rho)}| H |^{2}dv=0,$$
where $f$ is the immersion, $B_\rho$ is a geodesic ball of $N$ with radius $\rho$  (G. Wheeler \cite{Wh}).

\item Complete biharmonic submanifolds (resp., hypersurfaces) $M$ in a Riemannian manifold of non-positive sectional (resp., Ricci) curvature whose mean curvature vector satisfies $\int_M | H^p |dv<\infty$ for some $p>0$  (Y. Luo \cite{Luo2}).

\item Complete biharmonic submanifolds (resp., hypersurfaces) in a Riemannian manifold
whose sectional curvature (resp., Ricci curvature) is non-positive with at most polynomial
volume growth (Y. Luo \cite{Luo2}).

\item Complete biharmonic submanifolds (resp., hypersurfaces) in a negatively curved Riemannian
manifold whose sectional curvature (resp., Ricci curvature) is smaller that $-\epsilon$ for some
$\epsilon>0$ (Y. Luo \cite{Luo2}).

\item Proper $\epsilon$-superharmonic submanifolds $M$ with $\epsilon>0$ in a complete Riemannian manifold $N$ whose mean curvature vector satisfying the growth condition
$$\lim_{\rho\to \infty}\frac{1}{\rho^2}\int_{f^{-1}(B_\rho)}| H |^{2+a}dv=0,$$
where $f$ is the immersion, $B_\rho$ is a geodesic ball of $N$ with radius $\rho$, and $a\geq 0$ (Luo \cite{Luo2}).

\end{itemize}

On the other hand, it was proved by Y.-L. Ou and L. Tang in \cite{OT}  that the generalized Chen's biharmonic conjecture is false in general by constructing foliations of proper biharmonic hyperplanes in a $5$-dimensional conformally flat space with negative sectional curvature.

Further counter-examples were constructed  in \cite{LO} by T. Liang and Y.-L. Ou.

\subsection{Two related biharmonic conjectures}

Now, I present two  biharmonic conjectures related to my original biharmonic conjecture.
\vskip.06in

\hskip.3in {\bf Biharmonic Conjecture for Hypersurfaces}:
\emph{Every biharmonic hypersurface of Euclidean spaces is minimal.}
\vskip.06in

The global version of my original biharmonic conjecture can be found, for instance, in \cite{AM,M13}.
\vskip.06in

\hskip.3in {\bf Global Version of Chen's biharmonic Conjecture}:
\emph{Every complete biharmonic submanifold of a Euclidean space is minimal.}
 
 \vskip.1in

\section{Null 2-type submanifolds.}

From the definition of null 2-type submanifolds and formula (1.3), it follows that the mean curvature vector $H$ of a null  2-type submanifold   satisfies the following simple  condition:
 \begin{align}\Delta H=cH,\end{align}
 for some non-zero constant $c$.   In fact beside biharmonic submanifolds, null 2-type
submanifolds and 1-type submanifolds are the only submanifolds of a Euclidean space whose mean curvature vector satisfies condition (5.3) for some constant $c$ (Lemma 1 of \cite{C14}). Hence, null 2-type submanifolds  (together with biharmonic submanifolds, minimal submanifolds of Euclidean space and minimal submanifolds
of hypersphere) are the  special class of  submanifolds  
which can be characterized by the simple geometric condition (5.3). Since the author had classified biharmonic surfaces and null 2-type surfaces in $\mathbb E^3$, the classification of surfaces in $\mathbb E^3$,  satisfying condition (5.3) had been done in the 1980s.   

   From the classification of finite type planar curves, we know that there exist no null 2-type curves in a plane. Furthermore  it was  proved in \cite{C14} that null 2-type curves in
Euclidean spaces are circular helices in a Euclidean 3-space with nonzero torsion.  For null 2-type surfaces, first, we know in \cite{C13} that circular cylinders are the only null 2 type surfaces in Euclidean 3-space. Moreover,  null 2-type surfaces in a Euclidean 4-space are helical cylinders if they have constant mean curvature \cite{C14}.  
The author doesn't know whether every null 2-type surface in a Euclidean
4-space has constant mean curvature.  Also it is easy to see that  formula (1.5)
implies that {\it hyperplanes, hyperspheres, and
hypersurfaces of null 2-type are the only hypersurfaces of
a Euclidean space which have constant mean curvature and
constant scalar curvature.\/}

 As a generalization of \cite{C13}, \cite{FL} used the same
method of \cite{C13} to study null 2-type hypersurfaces with
at most two distinct principal curvatures, 
 (see also \cite{FGL2} for null 2-type conformally flat
hypersurfaces of dimension $\not= 3$). 

For null 2-type submanifolds of codimension
2, the author proposed in \cite{C18} the following 

{\bf Problem 11.} {\it Study and classify null 2-type
submanifolds. In particular, classify all null 2-type
surfaces in 4-dimensional Euclidean space and in
4-dimensional pseudo-Euclidean spaces.}

\vskip.2in
\section{Finite Type Submanifolds in Homogeneous Spaces} \vskip.05in

Let $N$ be a compact connected Riemannian homogeneous manifold with irreducible iso\-tropy action. Let $G$ be the identity component of the group of all isometries of $N$. Then $G$ is a compact Lie group which acts on $N$ transitively. For each positive eigenvalue $\lambda$ of the Laplacian operator $\Delta$ we denote by
$m_{\lambda}$ the multiplicity of $\lambda$. Let $\phi_{1},\ldots,\phi_{m_{\lambda}}$ be an orthonormal
basis of the eigenspace $V_{\lambda}$ with eigenvalue $\lambda$. We define a map $x: M \rightarrow \mathbb E^{m_{\lambda}}$ by
\begin{align}x(p)=c(\phi_{1}(p),\ldots, \phi_{m_{\lambda}}(p)),\end{align}
where $c$ is a positive constant. 

For a suitable positive constant $c$, $x$ defines an isometric 1-type immersion of $N$ into $\mathbb E^{m_{\lambda}}$. If $\lambda$ is the $i$-th positive eigenvalue of $\Delta$,
then the  immersion $x$ is called the {\it i-th standard immersion\/} of $N$. For example, if $N$ is a standard $n$-sphere,  the first standard immersion is in fact the standard imbedding of the $n$-sphere as a standard hypersphere of $\mathbb E^{n+1}$. 

Similar to the studies of finite type submanifolds of a hypersphere of Euclidean space,  one may consider the following 

{\bf Problem 12.} {\it Let $N$ be an irreducible compact homogeneous manifold immersed in a Euclidean space $\mathbb E^N$ by its first standard immersion $\phi$ and $M$ a
submanifold of $N$. When $M$ is of finite type in $\mathbb E^N$ via $\phi$? In particular, when $M$ is of 1- or 2-type in $\mathbb E^N$ via $\phi$?}

 When  $N$ is a projective space $FP^m$ over a field $F=R, C\, or \, H$ equipped with a standard Riemannian metric,  this problem has been studied in \cite{BC2,BU,C4,C5,GR,MR,Ros1,Ros2,Ros3,UD1,Di1,Di2,Di3} among others.  See \cite{BN} for $N$ to be either the real Grassmannian $G^{R}(p,q)$ or the space $U(n)/O(n)$. More recent results were obtained in \cite{Di4,Di5,Di6}.

\end{document}